\documentclass{amsart}

\newtheorem{theorem}{Theorem}[section]

\newtheorem{lemma}[theorem]{Lemma}
\newtheorem{proposition}[theorem]{Proposition}
\newtheorem{corollary}[theorem]{Corollary}

\theoremstyle{definition}

\theoremstyle{remark}

\renewcommand{\labelenumi}{(\roman{enumi})}

\DeclareFontFamily{U}{wncy}{}
\DeclareFontShape{U}{wncy}{m}{n}{<->wncyr10}{}
\DeclareSymbolFont{mcy}{U}{wncy}{m}{n}
\DeclareMathSymbol{\Sh}{\mathord}{mcy}{"58}

\newcommand\mynote[1]{\marginpar{\ \\ \small \tt #1}}
\newcommand\bel[1]{{\mynote{#1}}\begin{equation}\label{#1}}
\newcommand\mylabel[1]{\label{#1}\mynote{#1}}

\newcommand{\ZZ}{\mathbb{Z}}
\newcommand{\QQ}{\mathbb{Q}}

\newcommand{\FF}{\mathbb{F}}

\newcommand{\GG}{\mathbb{G}}

\newcommand  {\shC}     {\mathcal{C}}

\newcommand  {\shF}     {\mathcal{F}}

\newcommand  {\shR}     {\mathcal{R}}

\newcommand  {\shP}     {\mathcal{P}}

\newcommand  {\foX}     {\mathfrak{X}}

\newcommand  {\cid}     {\mathfrak{c}}
\newcommand  {\coker}   {\operatorname{coker}}
\newcommand  {\Coh}     {\operatorname{Coh}}
\renewcommand{\cong}    {\equiv}

\newcommand  {\End}     {\operatorname{End}}

\newcommand  {\Ext}     {\operatorname{Ext}}

\newcommand  {\GL}      {\operatorname{GL}}

\newcommand  {\id}      {\operatorname{id}}

\newcommand  {\dirlim}  {\varinjlim}
\newcommand  {\invlim}  {\varprojlim}

\newcommand  {\lra}     {\longrightarrow}

\newcommand  {\Mat}     {\operatorname{Mat}}
\newcommand  {\maxid}   {\mathfrak{m}}

\renewcommand{\O}       {\mathcal{O}}

\newcommand  {\Pic}     {\operatorname{Pic}}

\newcommand  {\pr}      {\operatorname{pr}}

\newcommand  {\ra}      {\rightarrow}

\newcommand  {\Spec}    {\operatorname{Spec}}
\newcommand  {\Spf}     {\operatorname{Spf}}

\newcommand  {\Sym}     {\operatorname{Sym}}

\def\mydate{\number\day\space\ifcase\month \or January\or February\or March\or 
April\or May\or June\or July\or
August\or September\or October\or November\or December\fi \space\number\year}

\DeclareFontFamily{U}{wncy}{}
\DeclareFontShape{U}{wncy}{m}{n}{<->wncyr10}{}
\DeclareSymbolFont{mcy}{U}{wncy}{m}{n}
\DeclareMathSymbol{\Sh}{\mathord}{mcy}{"58}

\usepackage{amscd,amssymb,amsmath}

\begin{document}

\title[Unipotent vector bundles on elliptic curves]
      {On the ring of unipotent vector bundles on elliptic curves in positive characteristics}

\author[Stefan Schr\"oer]{Stefan Schr\"oer}
\address{Mathematisches Institut, Heinrich-Heine-Universit\"at,
40225 D\"usseldorf, Germany}
\curraddr{}
\email{schroeer@math.uni-duesseldorf.de}

\subjclass[2000]{14H60, 14L05}

\dedicatory{26 August 2009}

\begin{abstract}
Using   Fourier--Mukai transformations, we prove some results about the ring of unipotent vector bundles
on   elliptic curves in positive characteristics.
This ring was determined by Atiyah in characteristic zero,
who showed that it is  a polynomial ring in one
variable. It turns out that the situation in characteristic $p>0$ is
completely different and rather bizarre:
the ring is nonnoetherian and contains a subring whose spectrum contains infinitely many
copies of $\Spec(\ZZ)$, which are glued  with successively higher and higher infinitesimal 
identification at the point corresponding to the prime $p$.
\end{abstract}

\maketitle
\tableofcontents
\renewcommand{\labelenumi}{(\roman{enumi})}

\section*{Introduction}
Using sheaf theoretic methods, 
Atiyah \cite{Atiyah 1957} described the category of locally free sheaves on 
elliptic curves. In particular, he determined
all \emph{irreducible} locally free sheaves.
Later, Oda  \cite{Oda 1971} gave another approach based on isogenies.
More recently, Polishchuk \cite{Polishchuk 2003} as well as Hein and Ploog \cite{Hein; Ploog 2005}
used  to this end the  Fourier--Mukai transformation, which was introduced by Mukai \cite{Mukai 1981}.
 
Given two irreducible sheaves $\shF,\shF'$ on an elliptic curve $E$, their tensor product is usually no longer irreducible,
and one has a decomposition $\shF\otimes\shF'=\bigoplus_i\shF_i^{\oplus \lambda_i}$
into irreducible summands.
To understand this decomposition, it esentially suffices to consider
sheaves $\shF$ having a filtration $0=\shF_0\subset\ldots\subset\shF_n=\shF$
with $\shF_i/\shF_{i-1}=\O_E$. Such sheaves are called \emph{unipotent}.
The classification of unipotent sheaves neither depends on the elliptic curve
nor the ground field:
For each integer $n\geq 1$, there is, up to isomorphism, precisely one irreducible
unipotent sheaf $\shF_n$ of rank $n$. For example, the extensions
of $\O_E$ by itself coming from the nontrivial classes in the 1-dimensional $H^1(E,\O_E)=\Ext^1(\O_E,\O_E)$
yield $\shF_2$. 

The decomposition of tensor products of unipotent sheaves, however, does depend on the characteristic.
Atiyah was able to determine the induced ring structure, in characteristic zero, 
of the free abelian group $R$ generated by the isomorphism classes of indecomposable unipotent sheaves:
It turns out that it is nothing but the polynomial ring $\ZZ[\shF_2]$.

The goal of this  paper is to analyze the ring $R$ in positive characteristics.
The idea is to use    \emph{Fourier--Mukai transformations},
which induce an equivalence between the category $\Coh_u(E)$ of unipotent sheaves and the 
category $\Coh_0(E)$ of   sheaves supported by the origin $0\in E$, according to
a result of Mukai \cite{Mukai 1978}.
Under this equivalence,   tensor products of unipotent sheaves corresponds to   
\emph{convolution products} of sheaves supported by the origin, which is defined in terms of the formal
group attached to the elliptic curve.
The latter has not much   to do with geometry, and is basically a matter of
linear and commutative algebra.
We shall see that this links our original problem to the    behavior of   Jordan normal forms
under tensor products, and to the modular representation ring
of the additive group of the  $p$-adic integers.

Our main results  are the following:
First, we show that the ring structure of  $R$ depends only on the characteristic $p\geq 0$,
and not on the elliptic curve. In particular,  there is no
difference between ordinary and supersingular elliptic curves.
Second, we show that $\Spec(R)$ contains infinitely many irreducible components,
such that $R $ is nonnoetherian.
Third, we describe the subring $R_\infty\subset R$ generated by 
the irreducible unipotent sheaves of prime power rank explicitly.
Its spectrum is an inverse limit of copies of $\Spec(\ZZ)$, which
are glued together at the points corresponding to the prime $p$
with respect to successively higher and higher infinitesimal identifications.
The localization $R_\infty\otimes\QQ$ is the 0-dimensional ring of almost
constant sequences in $(a_0,a_1,\ldots)$ of rational numbers, and its
spectrum is the Alexandroff compactification of a countable discrete space. 

The paper is organized as follows:
Section \ref{krull-schmidt fourier-mukai}  contains   basic facts
on the Krull--Schmidt Theorem and the Fourier--Mukai transformation.
In Section \ref{groups rings}, we introduce the convolution ring $R$ attached
to formal groups, and relate them to tensor products of unipotent sheaves and
the modular representation ring of the group of $p$-adic integers.
The basic properties of the multiplication table are gathered in Section \ref{multiplication table}.
I have included proofs in order to make the paper selfcontained, although the results are
probably not new.
The core of the paper is contained in Section \ref{application lucas}: Here we apply
a variant of Lucas' Theorem on congruences of binomial coefficients
to study the subring $R_\infty\subset R$  generated by elements with $p$-power rank
and give a detailed description of its spectrum.
In Section \ref{stone spaces} we interpret these results in terms of Stone spaces,
that is, compact totally discrete spaces.

%

\section{Krull--Schmidt and Fourier--Mukai}
\mylabel{krull-schmidt fourier-mukai}

We start by recalling   elementary facts about the Krull--Schmidt Theorem.
Let $\shC$ be an abelian category.
An object $V$ is called \emph{indecomposable} if $V\neq 0$, and in every decomposition $V=V_1\oplus V_2$
we have either $V_1=0$ or $V_2=0$.
If  every object admits a decomposition
$V=V_1\oplus\ldots\oplus V_n$ into indecomposable objects
and this decomposition is unique up to isomorphism and ordering of the summands, one says that the
\emph{Krull--Schmidt Theorem holds} for the category $\shC$.

According to Atiyah's result \cite{Atiyah 1956}, this is the case if each Hom set is endowed with the structure of
a finite-dimensional $k$-vector space so that the composition maps are bilinear.
A prominent example is the category of coherent sheaves with proper support
on a scheme or   algebraic stack of finite type over $k$, or on a
complex analytic space.

For lack of better notation, we define the  $K$-group
$K_0(\shC)$  as the free abelian group generated by the isomorphism classes
$[V]$ of objects in $\shC$, modulo the relations $[V]=[V']+[V'']$ coming from
direct sum decompositions $V\simeq V'\oplus V''$. Note that we use
only relations coming from split short exact sequences, rather than arbitrary short exact sequences.
If the Krull--Schmidt Theorem holds, then $K_0(\shC)$ is the free abelian
group generated by the isomorphism classes of indecomposable objects.
Often the category $\shC$ is endowed with tensor products, which induce
a ring structure on the abelian group $K_0(\shC)$. Given two indecomposable
objects $V',V''$, it is often a challenging problem to decompose $V'\otimes V''$
into indecomposable objects.

Let us now turn to abelian varieties.
Fix a ground field $k$ of characteristic $p\geq 0$, let $A$ be an abelian variety,
and consider the abelian category $\Coh(A)$ of coherent sheaves on $A$.
Inside, we have the subcategory $\Coh_0(A)$ of all coherent sheaves supported
by the origin $0\in A$, and the subcategory $\Coh_u(A)$ of all unipotent
sheaves. Here a coherent sheaf $\shF$ is called \emph{unipotent}
if there is a filtration $0=\shF_0\subset\shF_1\subset\ldots\subset\shF_n=\shF$
with subquotients $\shF_i/\shF_{i-1}\simeq\O_A$.
Note that unipotent sheaves are locally free,
and that tensor products of unipotent sheaves are unipotent.

A surprising result of Mukai   tells us that there is a canonical
equivalence of abelian categories $\Coh_0(A)\simeq\Coh_u(\hat{A})$, where $\hat{A}=\Pic^0_A$ is
the dual abelian variety.
Here it is most elegant to work with the triangulated category $D^b_\text{coh}(A)$ of 
bounded complexes of coherent sheaves on $A$.
Let $\shP$ be the Poincar\'e bundle on $A\times\hat{A}$.
The \emph{Fourier--Mukai transformation}   is defined as
$$
\Phi_\shP:D^b_\text{coh}(A)\lra D^b_\text{coh}(\hat{A}),\quad \shF^\bullet\longmapsto R\pr_{2*}(\shP\otimes\pr_1^*(\shF^\bullet)),
$$
where $\pr_i$ denote the projections from $A\times\hat{A}$.
According to Mukai's fundamental result (\cite{Mukai 1981}, Theorem 2.2), this functor is an equivalence of triangulated categories.
Up to shift by $g=\dim(A)$, it induces an equivalence of abelian categories
$$
\Coh_u(A)\lra\Coh_0(\hat{A}),\quad
\shF\longmapsto \Phi_\shP(\shF)[g] = R^g\pr_{2*}(\shP\otimes\pr_1^*(\shF)),
$$
as discussed in loc.\ cit.\ Example 2.9. Note that Mukai observed this equivalence already
in \cite{Mukai 1978}, Theorem 4.12. 
Moreover, he showed that the tensor product   corresponds under the Fourier--Mukai transformation 
to the convolution product, which
we shall discuss  in the next section.  For more on the Fourier--Mukai transformation,
we refer to Huybrecht's monograph \cite{Huybrechts 2006}.

The upshot is that to understand the decomposition of tensor products of unipotent sheaves into irreducible ones,
it suffices to understand the corresponding decomposition of convolution products for sheaves supported by
the origin. The latter is ``merely'' linear algebra, and in some aspects   a much simpler tasks.

\section{Formal groups and convolution rings}
\mylabel{groups rings}

Fix a ground field $k$ of characteristic $p\geq 0$, and let
$\foX=\Spf(k[[t]])$ be the formal affine line.
Consider the category $\Coh_0(\foX)$ of all coherent $\O_\foX$-modules $\shF$
supported by the origin $0\in \foX$. 
Then $V=\Gamma(\foX,\shF)$ is a finite dimensional vector space endowed
with a nilpotent endomorphism $\varphi:V\ra V$ induced by the action of the indeterminate $t\in k[[t]]$.
In fact, the functor
$$
\shF\longmapsto (V,\varphi)
$$
is an equivalence between   $\Coh_0(\foX)$ and
the category of finite-dimensional vector spaces  endowed
with a nilpotent endomorphism.
Using the Jordan normal form,
we conclude that for each indecomposable object is isomorphic to precisely
on $\shF_n$, $n\geq 1$, which corresponds to the $n$-dimensional vector space  $V=k[[t]]/(t^n)$
with the shift operator $\varphi(t^i)=t^{i+1}$, or equivalently the
standard vector space $V=k^n$ with the Jordan matrix
$$
J_n=\begin{pmatrix}
0\\
1 & 0 \\
  & \ddots & \ddots\\
  &        & 1      & 0
\end{pmatrix}
\in\Mat(n,k).
$$
Let $f_n=[\shF_n]$ be its isomorphism class, such that $K_0(\Coh_0(\foX))$ is the
free abelian group generated by the $f_n$, $n\geq 1$.
How to decompose  a given $\shF\in\Coh_0(\foX)$?
It turns out that the  Hilbert function 
$$
l_\shF(i)=\dim H^0(X,\shF/\maxid^{i}\shF)
$$
does the job:

\begin{lemma}
\mylabel{hilbert function}
Write $[\shF]=\sum\lambda_if_i$. Then we have $\lambda_i=2l_\shF(i)-l_\shF(i+1)-l_\shF(i-1)$.
\end{lemma}

\proof
By additivity, it suffices to treat the case $\shF=\shF_n$.  Set $l(i)=l_{\shF_n}(i)$.
We compute
$$
l(i)=\dim k[[t]]/(t^n,t^i)
=
\begin{cases}
i & \text{if $i\leq n$;}\\
n & \text{else.}
\end{cases}
$$
Whence the discrete derivative $l'(i)=l(i+1)-l(i)$ takes values
$$
l'(i)=
\begin{cases}
1 & \text{if $i\leq n-1$;}\\
0 & \text{else.}
\end{cases}
$$ 
In turn, the second discrete derivative $l''(i)=l(i+2)-2l(i+1)+l(i)$
is zero, except for $l''(n-1)=-1$. The assertion follows.
\qed

\medskip
Now choose a \emph{formal group law} $F(x,y)\in k[[x,y]]$, and regard $\foX$ as a \emph{formal group}.
Let $\mu:\foX\times\foX\ra\foX$ be the resulting multiplication morphism, 
which is given by $k[[t]]\ra k[[x,y]]$, $t\mapsto F$, and consider the \emph{convolution product}
$$
\shF\star\shF'=\mu_*(\shF\otimes_k\shF')
$$
for objects   $\shF\in\Coh_0(\foX)$. In terms of vector spaces with nilpotent endomorphisms,
the convolution product is given by $V\otimes V'$, with the induced nilpotent endomorphism
$F(\varphi\otimes\id_{V'},\id_V\otimes\varphi')$, which by abuse of notation we simply denote by $F(\varphi,\varphi')$.
It is easy to see that
the convolution product endows the free abelian group  $K_0(\Coh_0(\foX))$
with a ring structure.

Recall that by Lazard's results \cite{Lazard 1955},
the isomorphism classes of 1-dimensional formal groups over algebraically closed
ground fields of characteristic $p>0$ correspond to a single numerical invariant,
the \emph{height} $h\in\left\{\infty, 1,2,\ldots\right\}$.
The formal additive group $\hat{\GG}_a$, say with $F(x,y)=x+y$, has height $h=\infty$,
whereas the formal multiplicative group $\hat{\GG}_m$, say with $F(x,y)=x+y+xy$, has height $h=1$.
All other cases cannot be expressed in terms of  polynomials, and to my knowledge there  is no
\emph{explicit}   formula known for the formal group laws with $1<h<\infty$.
Somewhat surprising,  the formal group law   plays almost no role
for the multiplication table in our ring $K_0(\Coh_0(\foX))$:

\begin{proposition}
The multiplicative structure of the ring $K_0(\Coh_0(\foX))$
depends only on the characteristic $p\geq 0$ of our ground field $k$.
\end{proposition}

\proof
Set $\shF=\shF_m\star\shF_n$.
In light of Lemma \ref{hilbert function}, it suffices to check
that the Hilbert function $l_\shF(i)=\dim k[[x,y]]/(x^m,y^n,F^i)$ depends only on $p\geq 0$
rather then of $F$.
By definition, a formal group law $F(x,y)$ satisfies $F(x,0)=x$ and $F(0,y)=y$,
whence is of the form $F(x,y)=x+y+xyv=ux+y$ for some $u,v\in k[[x,y]]$
with $u$ invertible. Setting $x'=ux$, we  have
$$
k[[x,y]]/(x^m,y^n,F^i)= k[[x',y]/(x'^m,y^m,(x'+y)^i).
$$
The latter is clearly independent of the formal group law  $F$,
and in fact gives the Hilbert function with respect to the formal additive group $\hat{\GG}_a$.
Note that the dependence on the characteristic $p\geq 0$ enters
via the binomial expansion of $(x'+y)^i$ over $k$, as we shall see in the
next section.
\qed

\medskip
We call the ring $R=K_0(\Coh_0(\foX))$ the \emph{convolution ring} in characteristic $p\geq 0$.
Note that the underlying abelian group is free, with basis $f_i$, $i\geq 1$,
and that the multiplication table $f_mf_n=\sum\lambda_if_i$ depends only
on the characteristic $p\geq 0$. The unit element is $f_1=1$.
For some formulas, it is convenient to define $f_0=0$.
Clearly, the ring does not change under base field extensions $k\subset k'$.
The convolution ring is esentially the ring of unipotent locally free sheaves on elliptic curves
with respect to tensor products:

\begin{proposition}
The ring $K_0(\Coh_u(E))$ of unipotent
locally free sheaves on an elliptic curve $E$ in characteristic $p\geq 0$
is isomorphic to the convolution ring $R$ in characteristic $p$.
\end{proposition}

\proof
Let $\foX$ be the formal group attached to the elliptic curve $\hat{E}=E$.
Clearly, we have $\Coh_0(\foX)=\Coh_0(\hat{E})$.
According to \cite{Mukai 1978}, Theorem 4.12, the Fourier--Mukai transformation
$$
\Coh_u(E)\lra\Coh_0(\hat{E}),\quad
\shF\longmapsto \Phi_\shP(\shF)[g] = R^g\pr_{2*}(\shP\otimes\pr_1^*(\shF))
$$
is an equivalence of categories, and transforms tensor products into convolution products.
The assertion follows.
\qed

\medskip
In particular, it plays no role wether the elliptic curve $E$ is ordinary or supersingular.
This property, however, seems to enter if one considers the action of Schur functors $\Sym^p(\shF)$.

The multiplicities in $f_mf_n=\sum\lambda_if_i$
are related to a difficult problem in linear algebra:

\begin{proposition}
The multiplicities $\lambda_i$ are the number of Jordan blocks $J_i$
in the Jordan normal form of the endomorphism $J_m\otimes\id + \id\otimes J_n$ of $k^m\otimes k^n$.
\end{proposition}

\proof
We may choose the additive formal group law $F(x,y)=x+y$, and the assertion follows from
the definition of convolution products.
\qed

\medskip
This problem which was studied in characteristic zero, for example,
by Roth \cite{Roth 1934} and Trampus \cite{Trampus 1966}; for more recent developements in positive characteristics, see McFall \cite{McFall 1979} and 
Norman \cite{Norman 2008}.

\medskip
There is also an interesting connection to representation theory: Let $p>0$, and consider the additive profinite group $G=\ZZ_p=\invlim\ZZ/p^m\ZZ$
of $p$-adic integers. Let $\shR$ be the category of continuous representation of $G$
on finite-dimensional $\FF_p$-vector spaces. Such a representation $G\ra\GL(n,\FF_p)$ factors
over a finite quotient $\ZZ/p^m\ZZ$, and the image of the generator $1\in\ZZ/p^m\ZZ$
is a matrix $A\in\GL(n,\FF_p)$ whose minimal polynomial divides $T^{p^m}-1=(T-1)^{p^m}$.
Whence  the matrix $A-E_n$ is nilpotent, where $E_n\in\Mat(n,\FF_p)$
denotes the unit matrix. The upshot is that   $K_0(\shR)$ 
is the free abelian group generated by classes $f'_i$, $i\geq 1$
corresponding to the unipotent Jordan matrices $J_i+E_i$.
In turn, we  obtain a bijection of abelian groups
$$
R\lra K_0(\shR),\quad f_i\longmapsto f'_i.
$$
The tensor product endows $K_0(\shR)$ with a ring strucutre; we call it
the \emph{ring of continuous modular representations}.

\begin{proposition}
The preceding bijection respects multiplication,
and yields an identification of the convolution ring $R$ in characteristic $p$
with the  ring $K_0(\shR)$.
\end{proposition}

\proof
Here we use the multiplicative formal group law $F(x,y)=x+y+xy=(x+1)(y+1) -1$, together with
the fact that tensor products of unipotent matrices are unipotent.
Clearly, the number of unipotent Jordan blocks $J_i+E_i$ in the Jordan normal form
of the unipotent endomorphism $(J_m+E_m)\otimes(J_n+E_m)$ is the same as the number
of nilpotent Jordan blocks $J_i$ in the corresponding nilpotent endomorphism $F(J_m,J_n)=(J_m+E_m)\otimes(J_n+E_m)-E_{mn}$.
\qed

\medskip
This  ring of continuous modular representations was studied, for example,  
by Srinivasan \cite{Srinivasan 1964}, Green \cite{Green 1962}, and Ralley \cite{Ralley 1969}.

\section{The multiplication table}
\mylabel{multiplication table}

We keep the notation from the preceding section, and start to investigate
the multiplication table of the convolution ring $R$ in characteristic $p\geq 0$, whose underlying abelian group is
defined in terms of the formal affine line $\foX=\Spf(k[[t]])$.
Its multiplication is defined with the help of a formal group law $F(x,y)$,
but depends only on the characteristic $p\geq 0$. We therefore
assume from now on that the formal group law is simply $F(x,y)=x+y$.

Recall that $R$ is the free abelian group on the generators $f_n$, $n\geq 1$.
These classes are given by the coherent sheaves $\shF_n$, which in turn correspond
to the vector spaces $k[[t]]/(t^n)$ endowed with the shift operator,
or equivalently the vector space $k^n$ endowed with the Jordan matrix $J_n$.
We start with two obvious facts:

\begin{proposition}
\mylabel{dimension constraint}
Let $m,n\geq 0$ be integers, and write $f_mf_n=\sum\lambda_if_i$. Then  
$\sum i\lambda_i=mn$
\end{proposition}

\proof
Set $V=k[[t]]/(t^m)$ and $V'=k[[t]]/(t^n)$. Clearly $\sum i\lambda_i$ must be the dimension
of $V\otimes V'$, which indeed is $mn$.
\qed

\begin{proposition}
\mylabel{summand constraint}
Let $0\leq m\leq n$ be integers, and write $f_mf_n=\sum\lambda_if_i$. Then we have $\sum\lambda_i=m$.
\end{proposition}

\proof
The integer $\sum\lambda_i$ is the number of Jordan blocks in the Jordan normal form
of $F(J_m,J_n)$. Equivalently, it is the dimension of   $\shF/t\shF$, where $\shF=\shF_m\star\shF_n$.
The latter equals the dimension of 
$$
k[[x,y]]/(x^m,y^n,F)=k[[x,y]]/(x^m,y^n,x+y)
$$
which clearly equals $m$.
\qed

\medskip
The next observation permits us to translate certain problems from linear algebra into commutative algebra:

\begin{lemma}
\mylabel{largest smallest}
Let $m,n\geq 1$ be integers, and write $f_mf_n=\sum\lambda_if_i$.
Then the largest integer $r\geq 1$ with $\lambda_r\neq 0$ coincides with
the smallest integer $s\geq 1$ with the property $(x+y)^s\cong 0$ modulo $(x^m,y^n)$.
\end{lemma}

\proof
The largest integer $r\geq 1$ with $\lambda_r\neq 0$ is nothing but the size
of the largest Jordan block in the Jordan normal form of the nilpotent endomorphism $F(J_m,J_n)$
of $k^m\otimes k^n$. In turn, this is the smallest integer $s\geq 1$ with $F(J_m,J_n)^s=0$.
The linear map $k[[x,y]]/(x^m,y^n)\ra \End(k^m\otimes k^n)$ given by $x\mapsto J_m\otimes\id$
and $y\mapsto \id\otimes J_m$ is injective, which can be seen by identifying the vector space $k^m\otimes k^n$
with the $k$-algebra $k[[x,y]]/(x^m,y^n)$. Consequently, our $s\geq 1$ is also the smallest integer
with $F(x,y)^s\in (x^m,y^n)$.
\qed

\medskip
With this at hand we deduce several useful facts:

\begin{proposition}
\mylabel{nilpotence constraint}
Let $m,n\geq 0$ be integers, and write $f_mf_n=\sum\lambda_if_i$. Then  
$\lambda_i=0$ for all $i\geq m+n$.
\end{proposition}

\proof
By the Binomial Theorem,
$F(x,y)^{m+n-1}=(x+y)^{m+n-1}$ is contained in the ideal $(x^m,y^n)$,
and the assertion follows from Lemma \ref{largest smallest}.
\qed

\begin{proposition}
\mylabel{two relation}
We have
$$
f_2f_n=\begin{cases}
f_{n-1}+f_{n+1} & \text{if $p$ does not divide $n$;}\\
2f_n       & \text{else.}
\end{cases}
$$
for all $n\geq 1$.
\end{proposition}

\proof
The case $n=1$ is trivial, so assume $n\geq 2$.
According to Proposition \ref{summand constraint} and \ref{nilpotence constraint}, 
we have   $f_2f_n=f_i+f_j$ for some $i\leq j\leq n+1$, which furthermore
satisfies $i+j=2n$ by Proposition \ref{dimension constraint}.
The only solutions are $i=j=n$ and $i=n-1$, $j=n+1$.
Clearly  $(x+y)^n\cong\binom{n}{1}xy^{n-1}$ modulo $(x^2,y^n)$,
and therefore $(x+y)^n\cong 0$ if and only if $p$ divides $n$. Now the assertion
follows from Lemma \ref{largest smallest}.
\qed

\begin{corollary}
The ring $R$ is integral if and only if $p=0$. In this case, 
the homomorphism of rings $\ZZ[X]\ra R$ defined by $X\mapsto f_2$ is bijective.
\end{corollary}

\proof
If $p\neq 0$, then $(f_2-2f_1)f_p=0$, whence the ring $R$ is not integral.
Now suppose $p=0$. Using the relation $f_2f_n=f_{n-1}+f_{n+1}$, we inductively infer that
$f_2^n=f_{n+1}+\sum_{i\leq n}\lambda_if_i$ for some coefficients $\lambda_i$.
Consider the standard $\ZZ$-bases $t^j\in\ZZ[X]$ and $f_{j+1}\in R$, $j\geq 0$.
With respect to these bases, the matrix of $\ZZ[X]\ra R$, 
viewed as a homomorphism of free $\ZZ$-modules,
is an upper triangular matrix all whose diagonal entries are $1$, whence
is invertible.
\qed

\begin{corollary}
\mylabel{famous relation}
Let $1\leq m\leq n$ be integers. If $p=0$ or $p>m+n-2$,
then
$$
f_mf_n= \sum_{i=0}^{m-1} f_{m+n-1-2i}.
$$
\end{corollary}

\proof
By induction on $m\geq 1$. The case $m=1$ is trivial.
Now suppose $m\geq 2$, such that $f_{m}= f_2f_{m-1} -f_{m-2}$ by Proposition \ref{two relation}.
We have inductively
\begin{equation*}
\begin{split}
f_mf_n &=f_2f_{m-1}f_n-f_{m-2}f_n\\
       &=f_2\sum_{i=0}^{m-2}f_{m+n-2-2i} - \sum_{i=0}^{m-3}f_{m+n-3-2i}\\
       &=\sum_{i=0}^{m-2}(f_{m+n-3-2i} + f_{m+n-1-2i}) - \sum_{i=0}^{m-3}f_{m+n-3-2i}\\
       &=f_{n-m+1}+\sum_{i=0}^{m-2} f_{m+n-1-2i},\\
\end{split}
\end{equation*}
as desired.
\qed

\medskip
The following observation will be crucial in the next section:

\begin{proposition}
\mylabel{congruence constraint}
Let $m,n\geq 0$, and write $f_mf_n=\sum\lambda_if_i$.
Let $r\geq 0$ be an integer so that
$p$ divides $\binom{m+n-2-a}{m-1-b}$
for all $0\leq b\leq a\leq m+n-2-r$. Then $\lambda_i=0$ for all integers $i>r$.
\end{proposition}

\proof
We have
$$
(x+y)^r\cong \sum\binom{r}{j} x^jy^{r-j} \quad\text{modulo $(x^m,y^n)$},
$$
where the sum runs over all integers $j$ subject to the conditions
\begin{equation}
\label{conditions}
0\leq j\leq r,\quad j\leq m-1,\quad r-j\leq n-1.
\end{equation}
Setting  
$r=m+n-2-a$ and $j=m-1-b$,
we observe that the conditions in (\ref{conditions}) ensures that $0\leq b\leq a\leq m+n-2-r$.
Whence $(x+y)^r\cong 0$ modulo $(x^m,y^n)$.  Now Lemma \ref{largest smallest} implies that $\lambda_{r+1}=\lambda_{r+2}=\ldots=0$.
\qed

\section{Application of Lucas' Theorem}
\mylabel{application lucas}

Throughout this section, we assume that the characteristic is $p>0$.
We shall study product $f_mf_n$ in the convolution ring $R$ where the indices $m,n$ are
$p$-powers.
Our results hinge on the classical Theorem of Lucas   on congruences for binomial
coefficients (\cite{Lucas 1878}, Section 21).
Let $a,b\geq 0$ be integers, and $a=\sum a_ip^i$ and $b=\sum b_ip^i$ be
their $p$-adic expansion, with digits $0\leq a_i,b_i<p$. Then Lucas Theorem asserts
$$
\binom{a}{b} \cong \prod_i\binom{a_i}{b_i} \mod p.
$$, which we define
 For our purposes, the following variant is useful: Fix a prime power $q=p^i$,
and now consider the $q$-adic expansion $a=\sum a_iq^i$ and $b=\sum b_iq^i$,
with digits $0\leq a_i,b_i<q$.

\begin{lemma}
\mylabel{binomial congruence}
With the preceding notation, we also have 
$\binom{a}{b} \cong \prod_i\binom{a_i}{b_i}$ modulo $p$.
\end{lemma}

\proof
This can be seen as in the prove of Lucas' Theorem given by Fine \cite{Fine 1947}. 
We reproduce the argument for the sake of the reader: Fix $a=\sum a_iq^i$. Inside the polynomial ring $\FF_p[T]$, we have
$$
(1+T)^a=\prod_{i\geq 0}(1+T^{q^i})^{a_i}=\prod_{i\geq 0}\sum_{c\geq 0}\binom{a_i}{c}T^{cq^i} =
\sum_{c_0,c_1,\ldots\geq 0} \prod_{i\ge 0}\binom{a_i}{c_i}T^{c_iq^i}
$$
In the latter sum, a summand vanishes if   $c_i\geq a_i$ for  one index $i\geq 0$,
an in particular if $c_i\geq p$ for one $i\geq 0$.
Discarding these trivial summands and using the uniqueness of $q$-adic expansions  with digits $0\leq b_i<q$ for the numbers $b=\sum b_iq^i$, the above
expression equals
$\sum_{b\geq 0} \prod_{i\geq 0} \binom{a_i}{b_i}T^b$.
The result follows by comparing with the binomial expansion of $(1+T)^a$.
\qed

\medskip
The following congruence property of successive binomial coefficients will be crucial for us:

\begin{lemma}
\mylabel{lucas application}
Let $q=p^\nu$   and $0\leq m\leq q$ be integers. Then  
$\binom{m+q-2-a}{q-1-b}\cong 0$ modulo $p$ for all integers $0\leq b\leq a\leq m-2$.
\end{lemma}

\proof
Clearly, the first digit of the $q$-adic expansion of $m+q-2-a$
and $q-1-b$ are $m-2-a$ and $q-1-b$, respectively.
The assumptions yield the inequality  $m-2-a<q-1-b$,
which ensures $\binom{m-2-a}{q-1-b}=0$. The statement now follows
from Lemma \ref{binomial congruence}.
\qed

\begin{proposition}
\mylabel{prime products}
Let $q=p^\nu$ and $0\leq m\leq q$ be integers. 
Then $f_mf_q=mf_q$.
\end{proposition}

\proof
Write $f_mf_q=\sum\lambda_if_i$.
In light of Lemma \ref{lucas application}, we may apply 
Proposition \ref{congruence constraint} with $n=r=q$ and deduce
 that $\lambda_i=0$ for $i>q$. Whence we have
the inequality $i\lambda_i\leq q\lambda_i$ for all $i$.
We also know $\sum\lambda_i=m$ and $\sum i\lambda_i=mq$.
Whence 
$$
mq=\sum i\lambda_i\leq\sum q\lambda_i=qm.
$$
It follows that all our inequalities $i\lambda_i\leq q\lambda_i$
are actually equalities. The latter implies $\lambda_1=\ldots=\lambda_{q-1}=0$,
and finally $\lambda_q=m$.
\qed

\medskip
These relations have rather strange consequences. Consider the subset
$$
S=\left\{p^jf_{p^i}\mid i,j\geq 0\right\}\subset R.
$$
 This is a multiplicative
subset of   by Proposition \ref{prime products}, and $p$ becomes invertible in the
localization $S^{-1}R$.
We thus have a canonical map $\ZZ[p^{-1}]\ra S^{-1}R$. 

\begin{proposition}
\mylabel{localization bijective}
The canonical map  $\ZZ[p^{-1}]\ra S^{-1}R$ is bijective.
\end{proposition}

\proof
To check that the map is surjective, it suffices to verify that $f_m/1\in S^{-1}R$ is in its image.
Chose some $q=p^\nu$ with $m\leq q$. Then $f_mf_q=mf_q$, whence $f_m/1=m/1$ lies in the image.
The map is also injective: Suppose $m/p^\nu\in \ZZ[p^{-1}]$ maps to zero. Then $p^jf_{p^i}m=0$ for some $i,j\geq 0$,
whence $m=0$.
\qed

\medskip
Given an integers $\nu\geq 0$, we now consider the subgroup $R_\nu\subset R$ generated by
the $f_1,f_p,\ldots, f_{p^\nu}\in R$. According to Proposition \ref{prime products},
this is a subring, whence a finite flat $\ZZ$-algebra of rank $\nu+1$.

\begin{proposition}
The ring $R_\nu$ is reduced.
\end{proposition}

\proof
Let $x=\sum_{i\geq j}\mu_if_{p^i}$ be a nonzero element, say with $\mu_j\neq 0$.
We   compute $x^2=\mu_j^2p^j f_{p^j} + \sum_{i>j}\mu'_if_{p^i}$
for certain integers $\mu'_i$. It follows that $x^2\neq 0$, and this implies that
$R_\nu$ is reduced.
\qed

\medskip
We now seek to understand the geometry of the map $\Spec(R_\nu)\ra\Spec(\ZZ)$.
To this end, we first turn our attention to the fiber ring $R_\nu\otimes\FF_p$.

\begin{proposition}
The $\FF_p$-algebra $R_\nu\otimes\FF_p$ is isomorphic to the local Artin
ring $\FF_p[x_1,\ldots,x_\nu]/(x_1,\ldots,x_\nu)^2$.
\end{proposition}
 
\proof
Proposition \ref{prime products} implies that
$(f_{p^i}\otimes 1)( f_{p^j}\otimes 1 )=0$ for $1\leq i,j\leq \nu$. 
Whence we obtain a homomorphism
$$
\FF_p[x_1,\ldots,x_\nu]/(x_1,\ldots,x_\nu)^2\lra
R_\nu\otimes\FF_p,\quad
x_i\longmapsto f_{p^i}\otimes 1,
$$
which is obviously surjective.
This map is bijective, because  both algebras have the same dimension
as vector spaces over $\FF_p$.
\qed

\medskip
Next we look at the rings $R_\nu[1/p]=R_\nu\otimes\ZZ[1/p]$ obtained by
inverting $p$, which we   view as an algebra over $\ZZ[1/p]$. The elements 
$$
e_i=f_{p^i}\otimes 1/p^i\in R_\nu[p^{-1}],\qquad 0\leq i\leq\nu
$$ 
form a $\ZZ[1/p]$-basis and satisfy the relations
$e_ie_j=e_j$ for $0\leq i\leq j\leq \nu$. These relations imply that the linear
maps
$$
\varphi_j:R_\nu[1/p]\lra \ZZ[1/p],\quad e_i\longmapsto\begin{cases}
1 & \text{if $i\leq j$,}\\ 
0 & \text{else}
\end{cases}
$$
are homomorphisms of algebras. In turn, we obtain a homomorphism of algebras
$$
\Phi_\nu:R_\nu[1/p]\lra\prod_{i=0}^\nu\ZZ[1/p],\quad
x\longmapsto (\varphi_0(x),\varphi_1(x),\ldots,\varphi_\nu(x)).
$$

\begin{proposition}
\mylabel{localized structure}
The homomorphism $\Phi_\nu:R_\nu[1/p]\ra\prod_{i=0}^\nu\ZZ[1/p]$  is bijective.
\end{proposition}

\proof
Both  $\ZZ[1/p]$-modules in question are free of rank $\nu+1$, whence it suffices to check that $\Phi_\nu$ is surjective.
Clearly, the images
\begin{equation}
\label{phi map}
\Phi_\nu(e_i) = (\underbrace{0,\ldots, 0}_{\text{$i$ entries}},1,\ldots,1)\in\prod_{i=0}^\nu\ZZ[1/p],
\quad
0\leq i\leq\nu
\end{equation}
are module generators, whence the result.
\qed

\begin{theorem}
The affine scheme $\Spec(R)$ contains infinitely many irreducible components.
In particular, the ring $R$ is not noetherian.
\end{theorem}

\proof
Seeking a contradiction, we assume that the spectrum of $R$ 
has only finitely many irreducible components.
Then the spectrum $X$ of the localization $R\otimes\QQ$ also has
only finitely many irreducible components, say $X=X_1\cup\ldots\cup X_\nu$.
Now consider the subalgebra  $R_\nu\otimes\QQ\subset R\otimes\QQ$,
and let $f:X\ra Y_\nu$ be the induced morphism of affine scheme.
The space $Y_\nu=\Spec(R_\nu\otimes\QQ)$ is discrete and contains $\nu+1$ points $y_0,y_1,\ldots,y_\nu\in Y_\nu$.
Without loss of generality we may assume that $f(X_i)\neq y_0$ for all $i$,
that is, $y_0\not\in f(X)$. On the other hand, the map $f:X\ra Y_\nu$ is
dominant, because $R_\nu\ra R$ is injective, contradiction.
\qed

\medskip
We now may use Proposition \ref{localized structure} to determine the normalization of the  reduced 
ring $R_\nu$.
Consider the induced map
$R_\nu\hookrightarrow R_\nu[p^{-1}] \ra \prod_{i=0}^\nu\ZZ[1/p]$.
According to Formula (\ref{phi map}), we have
\begin{equation}
\label{integral maps}
\Phi_\nu(f_{p^i}) = (\underbrace{0,\ldots, 0}_{\text{$i$ entries }},p^i,\ldots,p^i),
\end{equation}
whence there is a factorization $\Phi_\nu:R_\nu\ra\prod_{i=0}^\nu\ZZ$. 

\begin{proposition}
The inclusion $R_\nu\subset\prod_{i=0}^\nu\ZZ$ is the normalization of the reduced ring $R_\nu$
in its total fraction ring.
\end{proposition}

\proof
The map in question becomes bijective after inverting $p$, and the ring $\prod_{i=0}^\nu\ZZ$ is normal.
We infer that the normalization is of the form $\prod_{i=0}^\nu R_i\subset\prod_{i=0}^\nu\ZZ$ for
some subrings $R_i\subset\ZZ$. Since the ring $\ZZ$ contains only one subring, namely itself,
the inclusion in question must be the normalization of $R_\nu$.
\qed

\begin{proposition}
The image of the inclusion map $\Phi_\nu:R_\nu\ra\prod_{i=0}^\nu\ZZ$  is the subring consisting
of elements of the form $(b_0,\ldots,b_\nu)$ with $b_j\cong b_{j-1}$ modulo $p^j$ for $1\leq j\leq \nu$.
\end{proposition}

\proof
Clearly, the entries $b_0=\ldots=b_{i-1}=0$ and $b_i=\ldots=b_\nu=p^i$ of $\Phi_\nu(f_{p^i})$ satisfy
the conditions. By additivity, these conditions hold for all images.
Conversely, suppose we have entries with $b_j\cong b_{j-1}$ modulo $p^j$ for $1\leq j\leq \nu$.
By induction on $\nu$, we have $\Phi_{\nu-1}(x)=(b_0,\ldots,b_{\nu-1})$ for some
$x=\sum_{i=0}^{\nu-1}\lambda_i f_{p^i}$. Then
$\Phi_\nu(x)=(b_0,\ldots,b_{\nu-1},b_{\nu-1})$ and 
$\Phi_\nu(x+\lambda_\nu f_{p^\nu}) = (b_0,\ldots,b_{\nu-1}, b_{\nu-1}+\lambda p^\nu)$.
Since $b_\nu\cong b_{\nu-1}$ modulo $p^\nu$, some $\lambda_\nu$  solves the equation $b_\nu=b_{\nu-1}+\lambda_\nu p^\nu$.
Whence $(b_0,\ldots,b_\nu)$ lies in the image.
\qed

\medskip
Now recall that the \emph{conductor ideal} $\cid\subset R_\nu$
for the birational inclusion $R_\nu\subset\prod_{i=0}^\nu\ZZ$ is defined as the annihilator ideal of $\coker(\Phi_\nu)=(\prod_{i=0}^\nu\ZZ)/R_\nu$.
This is the largest ideal in $R_\nu$ that is at the same time an ideal in the overring$\prod_{i=0}^\nu\ZZ$.

\begin{proposition}
The conductor ideal $\cid\subset\prod_{i=0}^\nu\ZZ$ is the principal ideal generated
by $(p,p^2,\ldots,p^\nu,p^\nu)$.
\end{proposition}

\proof
Let $(b_0,\ldots,b_n)\in\cid$. Then we have $b_ia_i\cong b_{i-1}a_{i-1}$ modulo $p^i$, for all
$a_i\in\ZZ$, $1\leq i\leq\nu$. Setting $a_i=0$ and $a_{i-1}=1$ we deduce $b_{i-1}\in p^i\ZZ$.
Moreover,  $a_\nu=1$ and $a_{\nu-1}=0$ yields $b_\nu\in p^\nu\ZZ$.
Whence the conductor algebra is contained in the principal ideal defined by $(p,p^2,\ldots,p^\nu,p^\nu)$.
Conversely, it is easy to check that the latter element lies in the conductor ideal.
\qed

\medskip
We thus have an identification of residue rings
$$
(\prod_{i=0}^\nu\ZZ)/\cid = (\prod_{i=0}^{\nu-1}\ZZ/p^{i+1}\ZZ) \times\ZZ/p^\nu\ZZ,
$$
and   $\Phi_\nu(f_{p^i})\in\prod_{i=0}^\nu\ZZ$, $1\leq i\leq \nu$ vanishes modulo $\cid$. In turn, we obtain
a commutative diagram
$$
\begin{CD}
\prod_{i=0}^\nu\ZZ  @>>> (\prod_{i=0}^\nu\ZZ)/\cid\\
@AAA @AAA\\
R_\nu @>>>\ZZ/p^{\nu}\ZZ
\end{CD}
$$
where the map on the right is the diagonal map, and the lower map is given by
sending $f_{p^i}$, $1\leq i\leq\nu$ to zero.
According to general properties of conductor ideals, this diagram is cartesian, and the induced commutative diagram of affine schemes
$$
\begin{CD}
\bigcup_{i=0}^\nu\Spec(\ZZ) @<<< \Spec((\prod_{i=0}^\nu\ZZ)/\cid)\\
@VVV @VVV\\
\Spec(R_\nu) @<<<\Spec(\ZZ/p^\nu\ZZ)
\end{CD}
$$
is cartesian and cocartesian \cite{Ferrand 2003}. Roughly speaking, the scheme $\Spec(R_\nu)$
is obtained from $\nu+1$ disjoint copies of $\Spec(\ZZ)$ by gluing the copies
along the points of characteristic $p$ together, with higher and higher infinitesimal
identification along   successive copies.

\section{Stone spaces}
\mylabel{stone spaces}

Keeping the notation from the previous section, we now consider the subring
$$
R_\infty=\bigcup_{\nu\geq 0} R_\nu\subset R,
$$
which is a free $\ZZ$-module generated by all the $f_{p^i}\in R$ with $i\geq 0$.
We seek to understand its spectrum.
To achieve this, we merely have to analyze the canonical inclusion $R_\nu\subset R_{\nu+1}$.
To simplify things, we shall  tensor with $\QQ$.
Consider the diagram
\begin{equation}
\label{inclusion maps}
\begin{CD}
R_\nu @>>> R_{\nu+1}\\
@V\Phi_\nu VV @VV\Phi_{\nu+1}V\\
\prod_{i=0}^\nu\QQ @>>> \prod_{i=0}^{\nu+1}\QQ
\end{CD}
\end{equation}
where the lower map is given by $(a_0,\ldots,a_\nu)\mapsto (a_0,\ldots,a_\nu,a_\nu)$, that is,
by duplicating  the last entry. It follows   from Formula (\ref{integral maps}) that this
diagram is commutative. Set $Y_\nu=\Spec(\prod_{i=0}^{\nu+1}\QQ)$.  The induced map 
$Y_{\nu+1}\ra Y_\nu$
is easy to understand: If $y^\nu_i\in Y_\nu$ denotes the point corresponding
to the $i$-th projection $\prod_{i=0}^\nu\QQ\ra\QQ$, then we have
$$
y^{\nu+1}_i\longmapsto\begin{cases}
y^\nu_i & \text{if $i\leq \nu$}\\
y^\nu_\nu& \text{if $i=\nu+1$}.
\end{cases}
$$
We now consider the inverse limit $Y=\invlim Y_\nu\subset\prod  Y_\nu$.
Being an inverse limit of finite discrete spaces,
it is a \emph{Stone space}, that is, a compact  and totally disconnected space.
In our case it  consists of the  points
$$
y_i=(y_0^0,y_1^1,\ldots,y_i^i,y^{i+1}_i,y_i^{i+2}\ldots)\in Y,\quad i\geq 0,
$$
whose entries become eventually constant in the lower index, together with the distinguished point
$$
y_\infty=(y_0^0,y_1^1,y_2^2,\ldots)\in Y,
$$
whose entries always change in the lower index. We now come back to our ring $R_\infty$:

\begin{theorem}
The      ring $R_\infty\otimes\QQ$ is 0-dimensional.  Its spectrum 
is homeomorphic, as a topological space, to $Y$, and this   is the Alexandroff compactification of the 
discrete space $\left\{y_0,y_1,\ldots\right\}$ obtained by putting the point  $y_\infty$ at ``infinity''.
\end{theorem}

\proof
Let $\pr_\nu:Y\ra Y_\nu$ be the canonical projection onto the $\nu$-th factor.
The subspace  $\left\{y_0,y_1,\ldots\right\}=Y\smallsetminus\left\{y_\infty\right\}$ is discrete,
because $\pr_{\nu}^{-1}(y_i^\nu)=\left\{y_i\right\}$ whenever $i<\nu$.
Since $Y$ is compact, it  must be the Alexandroff compactification of the
discrete subspace $Y\setminus\left\{y_\infty\right\}$,
by the  the uniqueness of the latter (\cite{TG 1-4}, \S9, No.\ 8, Theorem 4).

We clearly have $R_\infty=\dirlim_\nu R_\nu$, whence there is a canonical continuous map
$ \Spec(R_\infty)\ra \invlim_\nu\Spec(R_\nu)=Y$, and the latter is a homeomorphism by \cite{EGA IVc}, 
Corollary 8.2.10.
In particular, $R_\infty$ has Krull dimension zero.
\qed

\medskip
Reduced $0$-dimensional rings are also called \emph{absolutely flat},
or \emph{von Neumann regular}.
It is easy to give an explicit description of $R_\infty\otimes\QQ$.
Let us call a sequence of rational numbers $(a_0,a_1,\ldots)$ \emph{almost constant}
if $a_i=a_{i+1}=\ldots$ for some index $i\geq 0$. The description of $R_\infty$
given by (\ref{inclusion maps}) immediately gives:

\begin{proposition}
The ring $R_\infty\otimes\QQ$ is isomorphic to the subring of $\prod_{i=0}^\infty\QQ$
whose elements are the  almost constant sequences $(a_0,a_1,\ldots)$.
\end{proposition}

\medskip
Having understood the ring $R_\infty$, we turn our attention
to the ring extension $R_\infty\subset R$.
One might speculate wether this ring extension is integral.
At least, we have:

\begin{proposition}
The fiber of the morphism $\Spec(R)\ra\Spec(R_\infty)$
over the point $y_\infty$ consist of only one point.
\end{proposition}

\proof
For each $f_{p^i}\in R_\infty$, the open subset $D(f_{p^i})\subset\Spec(R_\infty)$
contains $y_\infty$ by (\ref{integral maps}), and we also have $y_\infty\in D(p)$.
The assertion now follows from Proposition \ref{localization bijective}.
\qed



\begin{thebibliography}{ccccc}

\bibitem{Atiyah 1956}
M.\ Atiyah:
On the Krull--Schmidt theorem with application to sheaves.
Bull.\ Soc.\ Math.\ France 84 (1956), 307--317. 

\bibitem{Atiyah 1957}
M.\ Atiyah:
Vector bundles over an elliptic curve.
Proc.\ London Math.\ Soc.\   7  (1957) 414--452. 

\bibitem{TG 1-4}
N.~Bourbaki:
Topologie g\'en\'erale. Chapitres 1--4. 
Hermann, Paris, 1971.

\bibitem{Ferrand 2003}
D.\ Ferrand:
Conducteur, descente et pincement.
Bull.\ Soc.\ Math.\ France  131  (2003), 553--585.

\bibitem{Fine 1947}
N.\ Fine:
Binomial coefficients modulo a prime.
Amer.\ Math.\ Monthly 54, (1947). 589--592. 

\bibitem{Green 1962}
J.\ Green:
The modular representation algebra of a finite group.
Illinois J.\ Math.\ 6 (1962) 607--619. 

\bibitem{EGA IVc}
A.\ Grothendieck:
\'El\'ements de g\'eom\'etrie alg\'ebrique IV: \'Etude locale des
sch\'emas et des morphismes de sch\'emas.
Publ.\ Math., Inst.\ Hautes \'Etud.\ Sci.\  28 (1966).

\bibitem{Hein; Ploog 2005}
G.\ Hein, D.\ Ploog:
Fourier--Mukai transforms and stable bundles on elliptic curves. 
Beitr\"age Algebra Geom.\ 46 (2005), 423--434. 

\bibitem{Huybrechts 2006}
D.\ Huybrechts:
Fourier--Mukai transforms in algebraic geometry.
Oxford Mathematical Monographs. The Clarendon Press, Oxford University Press, Oxford, 2006.

\bibitem{Lazard 1955}
M.\ Lazard:
Sur les groupes de Lie formels \`a un param\`etre.
Bull.\ Soc.\ Math.\ France 83 (1955), 251--274.

\bibitem{Lucas 1878}
E.\ Lucas:
Theorie des Fonctions Numeriques Simplement Periodiques. 
Amer.\ J.\ Math.\  1  (1878),  197--240.

\bibitem{McFall 1979}
J.\ McFall:
How to compute the elementary divisors of the tensor product of two matrices.
Linear and Multilinear Algebra 7 (1979), 193--201.

\bibitem{Mukai 1978}
S.\ Mukai:
Semi-homogeneous vector bundles on an abelian variety.
J.\ Math.\ Kyoto Univ.\ 18 (1978),  239--272. 

\bibitem{Mukai 1981}
S.\ Mukai:
Duality between $D(X)$ and $D(\hat X)$ with its application to Picard sheaves.
Nagoya Math.\ J.\ 81 (1981), 153-175.

\bibitem{Norman 2008}
C.\ Norman:
On Jordan bases for the tensor product and Kronecker sum and their elementary divisors over fields of prime characteristic.  
Linear Multilinear Algebra  56  (2008),   415--451.

\bibitem{Oda 1971}
T.\ Oda:
Vector bundles on an elliptic curve.
Nagoya Math.\ J.\ 43 (1971), 41--72. 

\bibitem{Polishchuk 2003}
A.\ Polishchuk:
Abelian varieties, theta functions and the Fourier transform.
Cambridge Tracts in Mathematics 153. 
Cambridge University Press, Cambridge, 2003.

\bibitem{Ralley 1969}
T.\ Ralley:
Decomposition of products of modular representations.  
J.\ London Math.\ Soc.\  44  (1969), 480--484.

\bibitem{Roth 1934}
W.\ Roth:
On direct product matrices.
Bull.\ Amer.\ Math.\ Soc.\ 40 (1934), 461--468. 

\bibitem{Srinivasan 1964}
B.\ Srinivasan:
The modular representation ring of a cyclic $p$-group.  
Proc.\ London Math.\ Soc.\   14  (1964), 677--688.

\bibitem{Trampus 1966}
A.\ Trampus:
A canonical basis for the matrix transformation $X\rightarrow AX+XB$.  
J.\ Math.\ Anal.\ Appl.\  14  (1966), 242--252.

\end{thebibliography}
\end{document}